\documentclass [leqno] {amsart}
  
\usepackage {amssymb}
\usepackage {amsthm}
\usepackage {amscd}
\usepackage{bm}

\usepackage{pstricks}
\usepackage{color}
\usepackage{graphicx}

\usepackage[all]{xy}
\usepackage {hyperref} 

\begin {document}

\bibliographystyle{alpha}
\theoremstyle{plain}
\newtheorem{proposition}[subsubsection]{Proposition}
\newtheorem{lemma}[subsubsection]{Lemma}
\newtheorem{corollary}[subsubsection]{Corollary}
\newtheorem{thm}[subsubsection]{Theorem}
\newtheorem{introthm}{Theorem}
\newtheorem*{thm*}{Theorem}

\theoremstyle{definition}
\newtheorem{definition}[subsubsection]{Definition}
\newtheorem{notation}[subsubsection]{Notation}

\newtheorem{example}[subsubsection]{Example}

\theoremstyle{remark}
\newtheorem{remark}[subsubsection]{Remark}

\numberwithin{equation}{subsubsection}


\newcommand{\arir}{\ar@{^{(}->}}
\newcommand{\aril}{\ar@{_{(}->}}
\newcommand{\are}{\ar@{>>}}

\newcommand{\xr}[1] {\xrightarrow{#1}}


\newcommand{\mc}[1]{\mathcal{#1}}


\newcommand{\Spec} {{\rm Spec}}



\newcommand{\Z} {\mathbb{Z}}
\newcommand{\Q} {\mathbb{Q}}

\newcommand{\F} {\mathbb{F}}

\newcommand{\OO}{\mathcal{O}}
\newcommand{\Fr}{{\rm Frob}}
\newcommand{\coker}{{\rm coker}}
\newcommand{\assumptionk}{Let $k$ be a perfect field of positive characteristic $p$. }
\newcommand{\perf}{{\rm perf}}

\title{On the Frobenius stable part of Witt vector cohomology}

\author{Andre Chatzistamatiou}
\address{Fachbereich Mathematik \\ Universit\"at Duisburg-Essen \\ 45117 Essen, Germany}
\email{a.chatzistamatiou@uni-due.de}

\thanks{This work has been supported by the SFB/TR 45 ``Periods, moduli spaces and arithmetic of algebraic varieties''}

\begin{abstract}
For a proper (not necessarily smooth) variety over a finite field with q elements, Berthelot-Bloch-Esnault proved
a trace formula which computes the number of rational points modulo q in terms of the Witt vector cohomology.
We show the analogous formula for Witt vector cohomology of finite length. In addition, we prove a vanishing
result for the compactly supported \'etale cohomology of a constant p-torsion sheaf on an affine Cohen-Macaulay 
variety. 
\end{abstract}

\maketitle

\section*{Introduction}
Let $p$ be a prime number and let $k$ be a finite field with $p^a$ elements. For a proper (not necessarily smooth) 
scheme $X$ over $k$ we know (at least) two congruence formulas for the number of $k$-rational points $\#X(k)$ modulo powers of $p$.
The first formula, by Katz \cite{Katz}, states that 
$$
\sum_{i\geq 0} (-1)^i {\rm Tr}(F^a\mid H^i(X,\OO_X)) \equiv \#X(k) \mod p.
$$ 
Here $F$ denotes the absolute Frobenius with its Frobenius linear operation on cohomology, but $F^a$ is linear. 
The second formula is due to Bloch-Illusie in the smooth case and Berthelot-Bloch-Esnault \cite{BBE} in general. 
If $W\OO_X=\varprojlim_n W_n\OO_X$ denotes
the sheaf of Witt vectors on $X$ with the Frobenius endomorphism $F$ then 
\begin{equation*}\label{introduction-equation-traceWitt}
\sum_{i\geq 0} (-1)^i {\rm Tr}(F^a\mid H^i(X,W\OO_X)\otimes \Q) \equiv \#X(k) \mod p^a.
\end{equation*}
In this paper we study this trace formula on a finite level, i.e. for $W_n\OO_X$ where $n\geq 1$ is an integer. For a fixed $n$ we obtain the following result. 
\begin{introthm}[cf. Corollary \ref{corollary-congruenceformula}] \label{introduction-thm-congruenceformula}
If the cohomology groups  $H^i(X,W_n\OO_X)$ are \emph{free} $W_n(k)$-modules for all $i\geq 0$ then 
\begin{equation*}\label{introduction-equation-congruenceformula}
\sum_{i\geq 0} (-1)^i {\rm Tr}(F^a\mid H^i(X,W_n\OO_X)) \equiv \#X(k) \mod p^{\min\{a,n\}}.
\end{equation*}
\end{introthm}
 However, the assumption of the theorem is non-trivial for $n\geq 2$. In particular, it implies that the Frobenius is bijective on $H^i(X,W_n\OO_X)$ for all $i$.  

More generally, the purpose of this paper is to study the maximal subspace of $H^*(X,W_n\OO_X)$ (and $H^*(X,W\OO_X)$) on which the Frobenius is a bijection,
we call it the \emph{Frobenius stable} Witt vector cohomology. By using compactifications we extend the definition of Frobenius stable Witt vector cohomology 
to separated schemes of finite type over a perfect field $k$ (${\rm char}(k)=p$), we denote these groups by $H^*_c(X,W_n\OO_X)_s$ and $H^*_c(X,W\OO_X)_s$
(Definition \ref{compcoh-definition}). 
In contrast to usual Witt vector cohomology (for proper schemes), the groups $H^*_c(X,W\OO_X)_s$ are always finitely generated $W(k)$-modules.  
The first result is a weak Lefschetz-type statement.
\begin{introthm}[cf. Theorem \ref{thm-weakLefschetz}] \label{intro-thm-weakLefschetz}
\assumptionk Let $X$ be an \emph{affine} scheme of finite type over $k$.
Suppose $X$ is equidimensional of dimension $d$ and suppose that $X$ is 
Cohen-Macaulay. Then 
$$
H^i_c(X,\OO_X)_s=0\quad \text{for all $i\neq d$.}
$$
\end{introthm}
The $W_n(k)$-modules $H^*_c(X,W_n\OO_X)_s$ together with the Frobenius endomorphism are closely related to compactly supported \'etale cohomology 
$H^*_{\text{\'et},c}(X\times_{k} \bar{k},\Z/p^n)$ equipped with the operation of the 
Galois group ${\rm Gal}(\bar{k}/k)$ (see Proposition \ref{proposition-competale} for a precise statement).    Via this correspondence the Theorem \ref{intro-thm-weakLefschetz} asserts that 
$$H^i_{\text{\'et},c}(X\times_{k} \bar{k},\Z/p)=0 \quad \text{for all $i\neq d$.}$$ 

We denote by $K$ the quotient field of $W(k)$. 
By using the comparison theorem of \cite{BBE} 
between Witt vector cohomology and compactly supported rigid cohomology
we prove that $H^*_c(X,W\OO_X)_s\otimes_{W(k)} K$ is the slope zero part of $H^*_{rig,c}(X/K)$. 

In order to prove Theorem \ref{introduction-thm-congruenceformula},
it is not sufficient to work with cohomology groups, but instead we have to work with perfect complexes. 
We observe that the Frobenius stable Witt vector cohomology of $X$ has a natural interpretation as Witt vector cohomology 
of the perfect scheme $X^{\perf}$ attached to $X$ (Proposition \ref{proposition-compperfect}). 
Suppose that $X$ is proper, 
then the \v{C}ech complex (associated to a finite affine covering) for $W_n(\OO_{X^{\perf}})$ is a perfect complex which we 
denote by $R\Gamma(X,W_n\OO_X)_s$. If the Frobenius acts bijectively on $H^*(X,\OO_X)$ 
then 
$$
{\rm Tr}(F^a\mid R\Gamma(X,W_n\OO_X)_s) = \sum_i (-1)^i {\rm Tr}(F^a\mid H^i(X,W\OO_X)\otimes \Q) \mod p^n.  
$$
This allows us to reduce Theorem \ref{introduction-thm-congruenceformula} to the trace formula of Berthelot-Bloch-Esnault \cite{BBE}.

\subsection*{Acknowledgements}
It is a pleasure to thank H\'el\`ene Esnault for her strong encouragement. 
This paper 
is inspired by the paper \cite{BBE}. In particular, the definition of Frobenius stable
Witt vector cohomology for non-proper varieties follows \cite{BBE}.

\section{The stable part of Witt vector cohomology}

\subsection{} \label{subsection-Dieudonnering}
Let $k$ be a perfect field of positive characteristic $p$. 
We denote by $W(k)$ the 
ring of Witt vectors and by $W_n(k)$ the ring of Witt vectors of length $n$.
The Frobenius automorphism on $W(k)$ and $W_n(k)$ is denoted by $\sigma$.
We denote by $\mathcal{D}$ the Dieudonn{\'e} ring
$
\mathcal{D}=W(k)[F,V]
$
with relations 
$$Fa=\sigma(a)F,\; aV=V\sigma(a),\; FV=VF=p,$$ 
for all $a\in W(k)$. For an integer $n\geq 1$, we define $\mc{D}_n:=\mc{D}/p^n$.

By a $\mathcal{D}$-module we will mean a left $\mathcal{D}$-module. 
Given a $\mc{D}$-module $M$ we obtain a $\sigma$-linear (and a $\sigma^{-1}$-linear) map $F:M\xr{} M$ (and $V:M\xr{} M$, respectively).

For a $W(k)$-module $M$
we define 
$$
\sigma_*M:= W(k)\otimes_{\sigma^{-1},W(k)} M.
$$ 
In other words, the multiplication in $\sigma_*(M)$ is twisted by $\sigma$: $a\otimes m=1\otimes \sigma(a)m$. If $M$ is equipped 
with a $\sigma$-linear map $F$ then $\sigma\otimes F$ defines a $\sigma$-linear map $\sigma_*M\xr{} \sigma_*M$. Similarly, if $M$ is a 
$\mc{D}$-module then $\sigma_*M$ inherits a $\mc{D}$-module structure.

\begin{definition}\label{definition-Frobstable}
Let $M$ be a $W(k)$-module together with a $\sigma$-linear map $F:M\xr{} M$. We define
$$
M_s:=\bigcap_{a\geq 1} F^a(M).
$$
We call $M_s$ the \emph{(Frobenius) stable} part of $M$. We call $M$ stable if $M_s=M$.
\end{definition}

Obviously, $F(M_s)\subset M_s$, and thus $M_s$ is equipped with the map $F$. If $M$ is a $\mc{D}$-module 
then $M_s$ is a $\mc{D}$-module.

\begin{proposition}\label{proposition-stablexact}
Let $M,M',M''$ be $W(k)$-modules equipped with $\sigma$-linear maps $F,F',F''$. 
Let $n\geq 1$ be an integer.
\begin{itemize}
\item[(i)] The equality $(\sigma_*(M))_s=\sigma_*(M_s)$ holds.
\item[(ii)] If $M$ is a finite and free $W_n(k)$-module then $M_s$ is a finite and 
free $W_n(k)$-module.
\item[(iii)] If $M$ is a finite $W(k)$-module then 
$
M_s\xr{\cong} \varprojlim_n (M/p^n)_s 
$
is an isomorphism.
\item[(iv)] If $M$ is a finite $W(k)$-module then $F:M_s\xr{} M_s$ is bijective. 
\item[(v)] Suppose that
$$
0\xr{} M' \xr{} M \xr{} M'' \xr{} 0
$$ 
is a short exact sequence (compatible with the $\sigma$-linear maps) and $M$ is a finite $W(k)$-module. Then 
$$
0\xr{} (M')_s \xr{} M_s \xr{} (M'')_s \xr{} 0
$$
is exact.
\end{itemize}
\begin{proof}
The assertion (i) is obvious. The statement (ii) follows from the elementary 
fact stated in Lemma \ref{lemma-stable/nil}. By assumption $k$ is perfect 
and this implies that $M_s$ and $M_{\rm nil}$ are sub-$W_n(k)$-modules. Therefore
$M_s$ is a projective $W_n(k)$-modules and thus free.   

For (iii): The $\sigma$-linear map on $M/p^n$ is induced by $M$. 
If $M$ is a finite $W(k)$-module then 
$M\xr{} \varprojlim_n M/p^n$ is an isomorphism. We get inclusions   
$$
M_s\subset  \varprojlim_n (M/p^n)_s \subset M.
$$
By Lemma \ref{lemma-stable/nil}, $F$ is bijective on $\varprojlim_n (M/p^n)_s$, 
thus $\varprojlim_n (M/p^n)_s \subset M_s$.

For (iv): Follows from (iii) and the fact that $F$ is bijective on $(M/p^n)_s$. 

For (v): We only need to prove that $M_s\xr{} (M'')_s$ is surjective. 
Consider the exact sequence 
$$
0\xr{} K_n \xr{} M/p^n \xr{} M''/p^n \xr{} 0.
$$
Lemma \ref{lemma-stable/nil} implies that 
$$
0\xr{} (K_n)_s \xr{} (M/p^n)_s \xr{} (M''/p^n)_s \xr{} 0
$$
is again exact. The projective system $(K_n)$ satisfies
the Mittag-Leffler condition, and taking the limit $\varprojlim$ implies the claim. 
\end{proof}
\end{proposition}

\begin{lemma}\label{lemma-stable/nil}
Let $M$ be an abelian group. Let $F:M\xr{} M$ be an endomorphism. 
Suppose that there are integers $n,m\geq 1$ such that 
$\ker(F^n)=\ker(F^{n+1})$ and $F^m(M)=F^{m+1}(M)$. Then 
there exists a unique decomposition 
$$
M=M_s\oplus M_{\rm nil}
$$  
with $F(M_s)\subset M_s, F(M_{\rm nil})\subset M_{\rm nil},$ such that the restriction of $F$ to
$M_s$ is bijective, and the restriction of $F$ to $M_{\rm nil}$ is nilpotent.
Moreover,
$$
M_s=F^{n+m}(M), \quad M_{\rm nil}=\ker(F^{n+m}).
$$


\end{lemma}

\begin{proof}
We leave the proof to the reader.
\end{proof}


\subsection{}
Let $X$ be a separated scheme of finite type over $k$, we denote by 
$W_n(\OO_X)$ the Witt sheaf of rings of $X$ of length $n$. We have the  
Frobenius endomorphism 
$$
F:W_n(\OO_X) \xr{} W_n(\OO_X), \quad (a_1,\dots,a_n)\mapsto (a_0^p,\dots,a_n^p).
$$
For $X=\Spec(k)$ we keep the notation $\sigma$ for $F$.  
Of course, $W_n(\OO_X)$ is an $W_n(k)$ module. We have the Verschiebung 
$$
V:W_n(\OO_X) \xr{} W_n(\OO_X), \quad (a_1,\dots,a_n)\mapsto (0,a_1,\dots,a_{n-1}),
$$
and the relation $V\circ F=F \circ V=p$.

If $X$ is proper then the cohomology groups $H^i(X,W_n(\OO_X))$ are finite $W_n(k)$ modules.
Moreover, $F$  (and $V$) induces a $\sigma$-linear map (a $\sigma^{-1}$-linear map, respectively,) 
$$
H^i(X,W_n(\OO_X))\xr{} H^i(X,W_n(\OO_X)), \quad \text{for all $i$.} 
$$
For simplicity we write $F=H^i(F)$ and $V=H^i(V)$.  

\subsection{} 
Let $X$ be a separated (not necessarily proper) scheme of finite type over $k$. 
We may choose a compactification $Y$ of $X$. 
Let $Z=Y\backslash X$ and choose a ideal sheaf $\mathcal{I}$ for $Z$. 
Denoting 
$$
W_n(\mathcal{I}):=\{(a_1,\dots,a_n)\in W_n(\OO_X); a_i\in \mathcal{I}\}
$$
we get by restriction a Frobenius and a Verschiebung endomorphism. 

The groups 
$H^i(Y,W_n(\mathcal{I}))$
are finite $W_n(k)$-modules and equipped with the $\sigma$-linear map $F$ and 
the $\sigma^{-1}$-linear map $V$ satisfying $FV=VF=p$.
Thus they are $\mc{D}_n$-modules.
 
In general, $H^*(Y,W_n(\mathcal{I}))$ depends on the choice of the ideal and the compactification.
However, we will see that the stable part does not.

\begin{lemma} \label{Wittindofcomp} Let $n\geq 1$ be an integer.
\begin{itemize}
\item[(i)] If $\mathcal{I}'$ is another ideal for $Z$ such that  
$\mathcal{I}'\subset \mathcal{I}$ then 
$$
H^i(Y,W_n(\mathcal{I}'))_s \xr{\cong} H^i(Y,W_n(\mathcal{I}))_s   
$$ 
is an isomorphism of $\mc{D}_n$-modules for all $i$.
\item[(ii)] Suppose $Y'$ is another compactification of $X$ with a morphism $g:Y'\xr{} Y$ which induces the identity on $X$. Then 
$$
g^*: H^i(Y,W_n(\mathcal{I}))_s \xr{} H^i(Y',W_n(\mathcal{I}\cdot \OO_{Y'}))_s
$$
is an isomorphism of $\mc{D}_n$-modules for all $i$.
\end{itemize}
\begin{proof}
There is a short exact sequence  
\begin{equation}\label{reductionsequence}
0 \xr{} W_{n-1}(\mathcal{I}) \xr{V} W_{n}(\mathcal{I}) \xr{R^{n-1}} \mathcal{I} \xr{} 0,
\end{equation}
where $V$ is the Verschiebung $(a_1,\dots,a_{n-1})\mapsto (0,a_1,\dots, a_{n-1})$ 
and $R^{n-1}$ is the projection $(a_1,\dots,a_{n})\mapsto a_1$. This gives a long exact sequence
of $\mc{D}_n$-modules: 
$$
\dots \xr{} H^{i-1}(Y,\mathcal{I}) \xr{} \sigma_*H^i(W_{n-1}(\mathcal{I})) \xr{V} H^i(W_{n}(\mathcal{I})) \xr{R^{n-1}} H^i(\mathcal{I}) \xr{} \dots
$$
The maps in (i) and (ii) are compatible with this long exact sequence, thus it is sufficient
to show the assertion in the case $n=1$. 

For (i). For $a\geq 1$ the morphism $\Fr^a: \mathcal{I}\xr{} \Fr^a_*  \mathcal{I}, r\mapsto r^{p^a},$ factors through 
\begin{equation}\label{Wittindofcompfactor}
\mathcal{I}\xr{\phi^a} \Fr^a_* (\mathcal{I}^{p^a}) \xr{\subset} \Fr^a_*  \mathcal{I},
\end{equation}
which induces a morphism of $\mc{D}_1$-modules ($V$ acts as zero)
$$
\phi^a:H^i(Y,\mathcal{I}) \xr{}  \sigma^a_* H^i(Y,\mathcal{I}^{p^a}).
$$
Since $F$ is bijective on the stable part,
$$
H^i(Y,\mathcal{I})_s \xr{\phi^a}  \sigma^a_* H^i(Y,\mathcal{I}^{p^a})_s 
\xr{\sigma^a_*(H^i({\rm inclusion}))_s}  \sigma^a_* H^i(Y,\mathcal{I})_s  
$$
is surjective. Since $\mathcal{I}^{p^a}\subset \mathcal{I}'$ for some $a$ this implies the 
surjectivity of
\begin{equation}\label{Wittindofcomp-map}
H^i(Y,\mathcal{I}')_s \xr{} H^i(Y,\mathcal{I})_s.
\end{equation}
We also obtain that $H^i(Y,\mathcal{I}^{p^a})_s\xr{} H^i(Y,\mathcal{I}')_s$ is surjective when $\mathcal{I}^{p^a}\subset \mathcal{I}'.$

In order to prove the injectivity of \ref{Wittindofcomp-map} it is sufficient to prove that $\phi^a$ is surjective.  
This follows from the commutative diagram 
$$
\xymatrix
{
H^i(Y,\mathcal{I})_s \ar[r]^{\phi^a} 
&
\sigma^a_* H^i(Y,\mathcal{I}^{p^a})_s 
\\
H^i(Y,\mathcal{I}^{p^a})_s \ar[u]^{H^i({\rm inclusion})_s} \ar[ur]_{F^a}
}
$$
and the surjectivity of $F^a$. 

For (ii). We prove the claim in two steps. In the first step we show that
\begin{equation}\label{Wittindofcomp1step}
H^i(Y,\mathcal{I})_s \xr{\cong} H^i(Y,g_*\mathcal{I}\OO_{Y'})_s,
\end{equation} 
and in the second step we prove
\begin{equation}\label{Wittindofcomp2step}
H^i(Y,g_*\mathcal{I}\OO_{Y'})_s \xr{\cong} H^i(Y',\mathcal{I}\OO_{Y'})_s.
\end{equation} 

Defining
$$
K_n=\ker(\mathcal{I}^n\xr{} g_*\mathcal{I}^n\OO_{Y'}), \quad 
C_n=\coker(\mathcal{I}^n\xr{} g_*\mathcal{I}^n\OO_{Y'}), 
$$
the $\oplus_{n\geq 0} \mathcal{I}^n$ modules $\oplus_{n \geq 0} K_n$ and 
$\oplus_{n\geq 0} C_n$ are finitely generated \cite[3.3.1]{EGA-IHES11-3.3.1}.
Thus there exists $d$ such that for all $n\geq d$ the following
maps are surjective
$$
\mathcal{I}\otimes_{\OO_Y} K_n \xr{} K_{n+1}, \quad \mathcal{I}\otimes_{\OO_Y} C_n \xr{} C_{n+1}. 
$$  
Since $K_n,C_n$ are supported in $Y\backslash X$  (for $n\leq d$), we can choose an integer $e\geq 1$ such that
$
\mathcal{I}^e K_n=0=\mathcal{I}^e C_n 
$
for all $n$.
It follows that for $m=e+d$ the maps 
\begin{equation}\label{Wittindcompvan1}
K_{n+m}\xr{} K_n, \quad C_{n+m}\xr{} C_n \quad \text{vanish for all $n$.}
\end{equation}

From the commutative diagram 
$$
\xymatrix
{
\Fr^a_* \mathcal{I} \ar[r]
&
\Fr^a_* g_*\mathcal{I}\OO_{Y'}
\\
\Fr^a_* \mathcal{I}^{p^a}\ar[u]^{\cup}\ar[r]
&
\Fr^a_* g_*\mathcal{I}^{p^a}\OO_{Y'} \ar[u]_{\cup}
\\
\mathcal{I} \ar@/^4pc/[uu]^{\Fr^a} \ar[r] \ar[u]^{\phi^a}
&
g_*\mathcal{I}\OO_{Y'} \ar@/_4pc/[uu]_{g_*\Fr^a} \ar[u]_{g_*\phi^a}
}
$$
we obtain induced morphisms $\Fr^a:K_1\xr{} \Fr^a_*K_1$ and 
$\Fr^a:C_1\xr{} \Fr^a_*C_1$ which factor through 
$\Fr^a_*(K_{p^a})\subset \Fr^a_*(K_{1})$ and $\Fr^a_*(C_{p^a})\xr{} \Fr^a_*(C_{1})$,
respectively.  Thus $\Fr$ is nilpotent on $K_1$ and $C_1$ by \ref{Wittindcompvan1}.
This proves \ref{Wittindofcomp1step}.

For \ref{Wittindofcomp2step} it is sufficient to show that $\Fr$ acts as a nilpotent endomorphism
on $R^ig_*(\mathcal{I}\OO_{Y'})$ for $i>0$. Then the assertion follows from the Leray spectral
sequence. In view of \cite[Appendix-Proposition~5]{Deligne-Appendix-Residues-duality} the map 
$
R^ig_*(\mathcal{I}^{p^a}\OO_{Y'}) \xr{} R^ig_*(\mathcal{I}\OO_{Y'}), 
$ 
induced by $\mathcal{I}^{p^a}\OO_{Y'} \subset \mathcal{I}\OO_{Y'}$, vanish if $a$ is sufficiently large. 
Thus the factorization \ref{Wittindofcompfactor} implies the claim.
\end{proof}
\end{lemma}

For two ideals of $Z=Y\backslash X$ we can take the intersection, and two compactifications 
can be dominated by a third one, thus Lemma \ref{Wittindofcomp} shows that 
$H^i(Y,W_n(\mathcal{I}))_s$
is independent of the choice of $\mathcal{I}$ and $Y$.

\begin{definition} \label{compcoh-definition}
Let $Y$ be a compactification of $X$ and $\mathcal{I}$ an ideal for $Y\backslash X$. 
For all $i$ we denote by $H^i_c(X,W_n\OO_X)_s$ the $\mc{D}_n$-module 
$H^i(Y,W_n(\mathcal{I}))_s$.
If $X$ is compact we will omit the index $c$.
\end{definition} 

\subsection{}
Let $X$ be of finite type and separated over $k$. Choose a compactification $Y$ and an 
ideal $\mathcal{I}$ for $Y\backslash X$. The restriction 
of the Frobenius to the nilradical $\mathcal{N}$ of $\OO_Y$ is nilpotent and from 
the short exact sequence 
$$
0 \xr{} \mathcal{I}\cap \mathcal{N} \xr{} \mathcal{I} \xr{} \mathcal{I}\OO_{Y_{{\rm red}}} \xr{} 0
$$ 
we conclude 
$$
H^i_c(X,W_n\OO_X)_s \xr{\cong} H^i_c(X_{{\rm red}},W_n(\OO_{X_{{\rm red}}}))_s \quad \text{for all $i,n$.} 
$$

Let $U\subset X$ be an open subset and choose an ideal $\mathcal{J}$ for $Y\backslash U$. In view of
the short exact sequence 
$$
0\xr{} \mathcal{J} \xr{} \mathcal{I} \xr{} \mathcal{I}/\mathcal{J} \xr{} 0 
$$
we get a long exact sequence
\begin{multline}\label{UXsequence}
\dots \xr{} H^i_c(U,W_n\OO_U)_s \xr{} H^i_c(X,W_n\OO_X)_s \xr{} H^i_c(X\backslash U,W_n\OO_{X\backslash U})_s \xr{} \\ H^{i+1}_c(U,W_n\OO_U)_s \xr{} \dots
\end{multline}

Let $U_1,U_2$ be open sets of $X$ such that $X=U_1\cup U_2$. Choose ideals $\mathcal{I}_1$ 
and $\mathcal{I}_2$ for $Y\backslash U_1$ and $Y\backslash U_2$, respectively. Then 
$\mathcal{I}_1+\mathcal{I}_2$ is an ideal for $Y\backslash X$ and $\mathcal{I}_1\cap \mathcal{I}_2$
is an ideal for $Y\backslash (U_1\cap U_2)$. From the short exact sequence 
$$
0\xr{} W_n(\mathcal{I}_1\cap \mathcal{I}_2) \xr{} W_n(\mathcal{I}_1)\oplus W_n(\mathcal{I}_2) \xr{} W_n(\mathcal{I}_1+ \mathcal{I}_2) \xr{} 0
$$
we get a long exact sequence
\begin{multline}\label{MV}
\dots \xr{} H^i_c(U_1\cap U_2,W_n\OO_{U_1\cap U_2})_s \xr{}  H^i_c(U_1,W_n\OO_{U_1})_s \oplus  H^i_c(U_2,W_n\OO_{U_2})_s \xr{}  \\
H^i_c(X,W_n\OO_{X})_s \xr{}  H^{i+1}_c(U_1\cap U_2,W_n\OO_{U_1\cap U_2})_s \xr{} \dots
\end{multline}

\begin{definition} 
Let $X$ be a separated scheme of finite type over a perfect field $k$ of positive characteristic $p$.
By taking the inverse limit we define
$$
H^i_c(X,W\OO_X)_s=\varprojlim_n H^i_c(X,W_n\OO_X)_s \quad \text{for all $i$.}
$$ 
If $X$ is compact we will omit the index $c$.
\end{definition}

We have a natural inclusion 
\begin{equation}\label{inclusionstable}
H^i_c(X,W\OO_X)_s=\varprojlim_n H^i(Y,W_n\mathcal{I})_s\subset \varprojlim_n H^i(Y,W_n\mathcal{I})=H^i(Y,W\mathcal{I}),
\end{equation}
and after inverting $p$ we obtain in the notation of \cite[2.13]{BBE}:
$$
H^i_c(X,W\OO_X)_s \otimes_{W(k)} K  \subset  H^i(Y,W\mathcal{I})\otimes_{W(k)} K =: H^i_c(X,W\OO_{X,K}).
$$
In general, the $\mc{D}$-module $H^i(Y,W\mathcal{I})$ is not a finite $W(k)$-module. However,
the $K$-vector space $H^i(Y,W\mathcal{I})\otimes_{W(k)} K$ is  finite dimensional and
independent of the choice of the compactification $Y$ and the ideal $\mathcal{I}$ (\cite[\textsection2]{BBE}).

\begin{proposition}
Vie the inclusion \ref{inclusionstable} we get 
$
H^i_c(X,W\OO_X)_s = H^i(Y,W\mathcal{I})_s  
$
for all $i$.
\begin{proof}
The inclusion $\supset$ follows immediately from the definitions, and $\subset$ follows since $F$ is bijective on $H^i_c(X,W\OO_X)_s$.
\end{proof}
\end{proposition}
 
\begin{proposition} \label{proposition-longexseq}
For all $n\geq 1$ there is a long exact sequence of $\mc{D}$-modules
\begin{equation*}
\begin{split}
\dots &\xr{} H^{i}_c(X,W\OO_X)_s \xr{p^n} H^{i}_c(X,W\OO_X)_s \xr{} H^{i}_c(X,W_n\OO_X)_s  \\
      &\xr{} H^{i+1}_c(X,W\OO_X)_s \xr{p^n} H^{i+1}_c(X,W\OO_X)_s \xr{} \dots.   
\end{split}
\end{equation*}
\begin{proof}
For all $m>n$ we have a short exact sequence (notation as in Definition \ref{compcoh-definition})
\begin{equation}\label{WmnIseq-longexseq}
0\xr{} W_{m-n}(\mathcal{I}) \xr{V^n} W_{m}(\mathcal{I}) \xr{R^{m-n}} W_n(\mathcal{I}) \xr{} 0,  
\end{equation}
with 
\begin{equation*}
\begin{split}
V^n(a_0,\dots,a_{m-n-1})&=(\underset{n}{\underbrace{0,\dots,0}},a_0,\dots,a_{m-n-1}), \\
R^{m-n}(a_0,\dots,a_{m-1})&=(a_0,\dots,a_{n-1}).
\end{split}
\end{equation*}
Now, \ref{WmnIseq-longexseq} induces a long exact sequence of $\mc{D}_m$-modules
$$
\dots \xr{} \sigma^n_* H^i(Y,W_{m-n}(\mathcal{I})) \xr{V^n} H^i(Y,W_{m}(\mathcal{I})) \xr{} H^i(Y,W_n(\mathcal{I})) \xr{} \dots.  
$$
Taking stable part and using the isomorphism 
$$
F^n:H^i(Y,W_{m-n}(\mathcal{I}))_s \xr{\cong}  \sigma^n_* H^i(Y,W_{m-n}(\mathcal{I}))_s, 
$$
we obtain a long exact sequence
\begin{equation}\label{equation-Wnlongexactseq}
\dots \xr{} H^i(Y,W_{m-n}(\mathcal{I}))_s \xr{V^n\circ F^n} H^i(Y,W_{m}(\mathcal{I}))_s \xr{} H^i(Y,W_n(\mathcal{I}))_s \xr{} \dots. 
\end{equation}
Since all groups are finite $W_m(k)$-modules the projective limit $\varprojlim_m$ is exact, and $V\cdot F=p$ implies the claim.
\end{proof}
\end{proposition}

\begin{corollary}\label{finitegeneration}
Let $X$ be a separated scheme of finite type over a perfect field $k$ of positive characteristic $p$.
The $W(k)$-module $H^i_c(X,W\OO_X)_s$ is finitely generated for all $i$.
\begin{proof}
Proposition \ref{proposition-longexseq} implies that 
$$
H^i_c(X,W\OO_X)_s/p^nH^i_c(X,W\OO_X)_s \subset H^i_c(X,W_n\OO_X)_s
$$
for all $i,n$. Therefore 
$$
H^i_c(X,W\OO_X)_s=\varprojlim_n \left( H^i_c(X,W\OO_X)_s/p^nH^i_c(X,W\OO_X)_s \right),
$$
and the assertion follows since $H^i_c(X,\OO_X)_s$ is a finite dimensional $k$-vector space. 
\end{proof}
\end{corollary}

\begin{remark}
Corollary \ref{finitegeneration} follows immediately from \cite[\textsection5, Proposition~3]{Serre}. We include a proof for the convenience of the reader.
\end{remark}

\begin{proposition}\label{proposition-FstableOimpliesFstableWO}
Suppose $X$ is proper over $k$. Suppose that $H^i(X,\OO_X)_s=H^i(X,\OO_X)$ for all $i$.
Then $H^i(X,W_n\OO_X)_s=H^i(X,W_n\OO_X)$ for all $i$ and $n\geq 1$. In particular, 
$H^i(X,W\OO_X)_s=H^i(X,W\OO_X)$ for all $i$.
\begin{proof}
The equality 
$ 
H^i(X,W_n\OO_X)_s=H^i(X,W_n\OO_X)
$
for all $i$, follows by induction on $n$ from the short exact sequence 
$$0 \xr{} W_{n-1}(\mathcal{O}_X) \xr{V} W_{n}(\mathcal{O}_X) \xr{R^{n-1}} \mathcal{O}_X \xr{} 0.$$
\end{proof}
\end{proposition}

\begin{proposition}\label{proposition-equivalencefreestable}
Let $X$ be a separated scheme of finite type over a perfect field $k$ of positive characteristic $p$. Let $n\geq 1$. The following statements are equivalent:
\begin{itemize}
\item[(i)] For all $i$, $H^i_c(X,W_n\OO_X)_s$ is a free $W_n(k)$-module. 
\item[(ii)] For all $i$, $H^i_c(X,W_n\OO_X)_s$ is a free $W_n(k)$-module of rank $\dim_k H^i_c(X,\OO_X)_s$. 
\item[(iii)] For all $i$, the map $H^i_c(X,W_n\OO_X)_s\xr{} H^i_c(X,\OO_X)_s$ is surjective.
\end{itemize}
\begin{proof}
Obviously (ii) implies (i). Now, suppose (i) holds. It is clear from the long exact sequence
\ref{equation-Wnlongexactseq} that 
\begin{equation}\label{equation-rankinequality}
{\rm length}(H^i_c(X,W_m\OO_X)_s)\leq m\cdot \dim_k H^i_c(X,\OO_X)_s
\end{equation}
for all $m\leq n$.
If the canonical map $H^i_c(X,W_n\OO_X)_s/p \xr{} H^i_c(X,\OO_X)_s$ is surjective 
(which holds for $i\geq \dim X$)
then equality holds in \ref{equation-rankinequality} for $m=n$, and thus
$$
H^i_c(X,W_{n-1}\OO_X)_s\xr{V\circ F} H^i_c(X,W_n\OO_X)_s
$$ 
is injective. It follows that $H^{i-1}_c(X,W_n\OO_X)_s/p \xr{} H^{i-1}_c(X,\OO_X)_s$ 
is surjective. By descending induction on $i$ (starting from $i=\dim X$) we see that (ii) holds.

Now, suppose that (iii) holds. By induction on $n$ we may suppose that for all $i$, 
$H_c^i(X,W_{n-1}\OO_X)_s$ is free of rank $\dim_k H_c^i(X,\OO_X)_s$. 
Since $H^i_c(X,W_{n-1}\OO_X)_s/p\xr{} H_c^i(X,\OO_X)_s$ is surjective (by assumption (iii))
it is an isomorphism. It follows that 
\begin{equation}\label{equation-WnWn-1surjective}
H^i_c(X,W_{n}\OO_X)_s\xr{} H^i_c(X,W_{n-1}\OO_X)_s
\end{equation}
is surjective. From the long exact sequence \ref{equation-Wnlongexactseq} we get short exact sequences
$$
0\xr{} H_c^i(X,W_{n-1}\OO_X)_s \xr{V\circ F} H_c^i(X,W_{n}\OO_X)_s \xr{} H^i_c(X,\OO_X)_s \xr{} 0.
$$ 
The surjectivity of \ref{equation-WnWn-1surjective} and $V\circ F=p$ on $H^i_c(X,W_{n}\OO_X)_s$ 
implies 
$$H_c^i(X,W_{n-1}\OO_X)_s=p\cdot H_c^i(X,W_{n}\OO_X)_s$$ 
which yields (together with $H^i_c(X,W_{n}\OO_X)_s/p\cong H^i_c(X,\OO_X)_s$) that $H^i_c(X,W_{n}\OO_X)_s$ is free of rank $\dim_k H^i_c(X,\OO_X)_s$. 
\end{proof}
\end{proposition}

\begin{proposition}(cf. \cite[\textsection5, Corollaire~2]{Serre})\label{proposition-properfree}
Suppose $X$ is proper and $n\geq 2$. The following statements are equivalent:  
\begin{itemize}
\item[(i)] For all $i$, $H^i(X,W_n\OO_X)$ is a free $W_n(k)$-module. 
\item[(ii)] For all $i$, $H^i(X,\OO_X)_s=H^i(X,\OO_X)$ and 
$H^i(X,W_n\OO_X)\xr{} H^i(X,\OO_X)$ is surjective.
\item[(iii)] For all $i$, $H^i(X,W_n\OO_X)_s=H^i(X,W_n\OO_X)$ and $H^i(X,W_n\OO_X)_s$ is a free $W_n(k)$-module.
\end{itemize}
\begin{proof}
Proposition \ref{proposition-FstableOimpliesFstableWO} and \ref{proposition-equivalencefreestable} imply
$(ii)\Rightarrow (iii)$. Obviously, $(iii)\Rightarrow (i)$. 

Now, suppose that (i) holds. From the long exact sequence associated to 
\ref{WmnIseq-longexseq} we obtain 
\begin{equation}\label{equation-rankinequalityW}
{\rm length}(H^i(X,W_m\OO_X))\leq m\cdot \dim_k H^i(X,\OO_X)
\end{equation}
for all $m$. If $H^i(X,W_n\OO_X)/p \xr{} H^i(X,\OO_X)$ is surjective then in \ref{equation-rankinequalityW} equality 
holds; this implies that 
$$
V:\sigma_*H^i(X,W_{n-1}\OO_X)\xr{} H^i(X,W_n\OO_X)
$$
is injective. Thus $V(\sigma_*H^i(X,W_{n-1}))=pH^i(X,W_n\OO_X)$ -- in particular, we see that $H^i(X,W_{n-1}\OO_X)$ is free. 
We also obtain $H^{i-1}(X,W_n\OO_X)/p\xr{\cong} H^{i-1}(X,\OO_X)$. By induction on $i$ we get 
$H^{i}(X,W_n\OO_X)/p\xr{\cong} H^{i}(X,\OO_X)$ for all $i$. We conclude that $H^i(X,W_m\OO_X)$ is free of rank 
$\dim_k H^i(X,\OO_X)$ for all $i,m\leq n$. Consider the case $m=2$: we have a short exact sequence
$$
0\xr{} \sigma_*H^i(X,\OO_X) \xr{V} H^i(X,W_2\OO_X) \xr{R} H^i(X,\OO_X) \xr{} 0 
$$ 
for all $i$. The composition $H^i(X,\OO_X) \xr{F} \sigma_*H^i(X,\OO_X)\xr{V} H^i(X,W_2\OO_X)$ equals $p\cdot R^{-1}$. 
Because $H^i(X,W_2\OO_X)$ is free, $V\circ F$ and thus $F$ is injective, which 
proves $H^i(X,\OO_X)_s=H^i(X,\OO_X)$.
\end{proof}
\end{proposition}

\begin{thm}\label{thm-weakLefschetz}
\assumptionk Let $X$ be an \emph{affine} scheme of finite type over $k$.
Suppose $X$ is equidimensional of dimension $d$ and suppose that $X$ is 
Cohen-Macaulay. Then 
$$
H^i_c(X,\OO_X)_s=0\quad \text{for all $i\neq d$.}
$$
\begin{proof}
Let $Y$ be a compactification of $X$. By blowing up the complement
$Y\backslash X$ we may suppose that we can find an ideal $\mc{I}$ for 
$Y\backslash X$ which is a Cartier divisor. It is sufficient to prove that 
the projective system 
$
(H^i(Y,\mc{I}^n))_n
$
is essentially zero if $i<d$. Let $\omega_Y$ be the dualizing complex of $Y$. We obtain
$$
H^i(Y,\mc{I}^n)^{\vee}=H^{-i}(Y,\omega_Y\otimes \mc{I}^{-n}),
$$
thus we need to prove that the inductive system $(H^{-i}(Y,\omega_Y\otimes \mc{I}^{-n}))_n$ is essentially
zero. Since 
$$
\varinjlim_{n} \omega_Y\otimes \mc{I}^{-n} = j_*({\omega_Y}_{\mid X}),
$$
with $j:X\xr{} Y$ the open immersion, we get 
\begin{equation*}
\begin{split}
\varinjlim_n H^{-i}(Y,\omega_Y\otimes \mc{I}^{-n}) &= H^{-i}(Y,j_*( {\omega_Y}_{\mid X})) \\
                                                   &= H^{-i}(Y,Rj_*({\omega_Y}_{\mid X})) \\
                                                   &= H^{-i}(X,{\omega_Y}_{\mid X}).
\end{split}
\end{equation*}
By assumption $X$ is Cohen-Macaulay, and therefore ${\omega_Y}_{\mid X} \cong \omega_X$ is concentrated in degree $-d$.
\end{proof}
\end{thm}

\begin{remark}
In general the statement of the Theorem \ref{thm-weakLefschetz} fails
if $X$ is not Cohen-Macaulay. For example consider two affine planes glued at a point $0$, 
$X=\mathbb{A}^2\cup_0 \mathbb{A}^2$. Then $H^1_c(X,\OO_X)_s\cong k$.
\end{remark}

\begin{corollary}
With the assumptions of Theorem \ref{thm-weakLefschetz}. The group $H^i_c(X,W\OO_X)_s$
vanishes if $i\neq d$, and is a finite free $W(k)$-module for $i=d$. 
\begin{proof}
Follows from Proposition \ref{proposition-longexseq} and \ref{proposition-equivalencefreestable}.
\end{proof}
\end{corollary}

\subsection{}
Next, we want to show that 
$$
H^i_c(X,W\OO_X)_s \otimes_{W(k)} K \cong H^i_{rig,c}(X/K)_{[0]} 
$$
where the right hand side is the slope $=0$ part of compactly supported rigid cohomology. 
We will need the following Lemma.

\begin{lemma}\label{lemmabasechange}
Suppose $k\subset L$ is an extension of perfect fields. 
Let $(M_n)$ be a projective system of $W(k)$-modules such that $p^nM_n=0$ and $M_n$ is a finite $W_n(k)$ module.
We set $M=\varprojlim_n M_n$ and consider the natural map 
\begin{equation}\label{modulebasechange}
M\otimes_{W(k)} W(L) \xr{} \varprojlim_n (M_n\otimes_{W_n(k)} W_n(L)).
\end{equation}
\begin{itemize}
\item[(i)] If $M$ is a finite $W(k)$-module then the map  \ref{modulebasechange} is an isomorphism.
\item[(ii)] The map \ref{modulebasechange} is injective.
\end{itemize}
\begin{proof}[Sketch of proof]
Since $k$ and $L$ are perfect we conclude that $W_n(L)$ is flat over $W_n(k)$.
Assertion (i) follows by reduction to the case $M_n=M/p^n$. Statement (ii) follows
from (i). We leave the details to the reader.
\end{proof}
\end{lemma}
 
\begin{proposition}\label{proposition-slopezero} 
Let $k$ be a perfect field of positive characteristic $p$, with ring of Witt vectors $W=W(k)$ and $K={\rm Frac}(W(k))$. 
Let $X$ be a separated scheme of finite type over $k$. 
For all $i$, there is an isomorphism 
$$
H^i_c(X,W\OO_X)_s \otimes_W K \xr{} H^i_{rig,c}(X/K)_{[0]}
$$
which is compatible with the $F$-operation.
\begin{proof}
Let $Y$ be a compactification of $X$ and $\mathcal{I}$ an ideal for $Y\backslash X$. 
By the work of Berthelot, Bloch and Esnault \cite{BBE}, there is an isomorphism
of $F$-isocrystals 
$$
H^i(Y,W(\mathcal{I}))\otimes_W K \xr{} H^i_{rig,c}(X/K)_{<1},
$$ 
where the right hand side is the slope $<1$ part of (compactly supported) rigid
cohomology. By definition we have
$$
H^i_c(X,W(\OO_X))_s \otimes_W K \subset H^i(Y,W(\mathcal{I}))\otimes_W K
$$
and since $F$ is bijective on the finite $W(k)$-module 
$H^i_c(X,W(\OO_X))_s$ we obtain 
\begin{equation}\label{maptoslopezero}
H^i_c(X,W(\OO_X))_s \otimes_W K \subset H^i_{rig,c}(X/K)_{[0]}.
\end{equation}

We set $\tilde{K}={\rm Frac}(W(\bar{k}))$. In order to prove  the surjectivity 
of \ref{maptoslopezero}, we need to show that 
every $v\in H^i_{rig,c}(X/K)_{[0]}\otimes_K \tilde{K}$ with 
$(F\otimes \sigma)(v)=v$ is contained in 
$H^i_c(X,W(\OO_X))_s \otimes_W \tilde{K}$.  
Multiplying $v$ by a suitable power of $p$ we may
assume that $v$ lies in the image of 
$
H^i(Y,W(\mathcal{I}))\otimes_W W(\bar{k}).
$
Choose a preimage $v'$ of $v$; $v'$ is well-defined up to $p$-power torsion.
Again, by multiplying $v'$ with a power of $p$ we may assume that 
$(F\otimes \sigma)(v')=v'$. Denote by $v'_n$ the image of $v'$ in 
$H^i(Y,W_n(\mathcal{I}))\otimes_{W_n} W_n(\bar{k})$, it is obviously
contained in the Frobenius stable part (for $F\otimes \sigma$). It is easy to see that
$$
(H^i(Y,W_n(\mathcal{I}))\otimes_{W_n} W_n(\bar{k}))_s = H^i(Y,W_n(\mathcal{I}))_s\otimes_{W_n} W_n(\bar{k}),
$$
and $v'\in H^i_c(X,W(\OO_X))_s \otimes_W W(\bar{k})$ follows from Lemma \ref{lemmabasechange}(i),(ii).
\end{proof}
\end{proposition}

\begin{corollary}
Under the assumption of Proposition \ref{proposition-slopezero}. If $H^i_c(X,\OO_X)_s=0$ for all $i$ then 
$H^i_{rig,c}(X/K)_{[0]}$ for all $i$.
\begin{proof}
The vanishing of $H^*_c(X,\OO_X)_s$ implies the vanishing of $H^*_c(X,W\OO_X)_s$ by Proposition \ref{proposition-longexseq} and
Corollary \ref{finitegeneration}.
\end{proof}
\end{corollary}

\section{Comparison with \'etale $p$-adic cohomology}

\subsection{} \assumptionk  We fix an algebraic closure $\bar{k}$. 
Let $G$ be the Galois group of $\bar{k}$ over $k$.
We denote 
by $\mc{D}(k)$ and $\mc{D}(\bar{k})$ the Dieudonn{\'e} ring of $k$ and $\bar{k}$, 
respectively (see Section \ref{subsection-Dieudonnering}).
We have a functor 
\begin{multline*} 
\mc{G}:(\text{Frobenius stable $\mc{D}(k)$-modules $M$ s.t. $M$ is a finite $W(k)$-module})
\\
\xr{} (\text{finite $\Z_p$-modules with $G$-action})
\end{multline*}
which is defined as follows. For $M$ we denote by 
$$
\bar{M}:=\mc{D}(\bar{k})\otimes_{\mc{D}(k)} M=W(\bar{k})\otimes_{W(k)} M
$$
the induced module over $\mc{D}(\bar{k})$. It is equipped with the
obvious $G$-action, and the $G$-action commutes with $F,V$. Moreover, we have $\bar{M}^G=M$. The functor 
$$
N\mapsto N^{1-F}:=\ker(1-F:N\xr{} N)
$$
is exact when restricted to the category of $\mc{D}(\bar{k})$-modules $N$ such that $N$ is a finite $W(\bar{k})$-module (see \cite[II, Lemme~5.3]{Illusie}).
Moreover, if $N$ is stable then 
\begin{equation} \label{backtoN}
W(\bar{k}) \otimes_{\Z_p} N^{1-F}  \xr{\cong} N,
\end{equation}
thus $N^{1-F}$ is a finite $\Z_p$-module ($W(\bar{k})$ is faithfully flat over $\Z_p$). We set 
$$
\mc{G}(M)=\bar{M}^{1-F}.
$$
The functor $\mc{G}$ is exact and fully faithful. 

Let $M$ be again a stable $\mc{D}(k)$-module such that $M$ is a finite $W(k)$-module. Suppose that 
$
M=\varprojlim_n M_n
$
for stable $\mc{D}_n(k)$-modules $M_n$ such that $M_n$ is a finite $W_n(k)$-module. Lemma \ref{lemmabasechange}(i)
implies that $\bar{M}=\varprojlim_n \bar{M}_n$ and thus 
\begin{equation}\label{compGprojlim}
\mc{G}(M)=\varprojlim_n \mc{G}(M_n).
\end{equation}


\begin{proposition}\label{proposition-competale}
Let $k$ be a perfect field of positive characteristic $p$ with ring of Witt vectors $W=W(k)$ and $K={\rm Frac}(W(k))$. Let $X$ be a separated scheme of finite type over $k$.
\begin{enumerate} 
\item For all $i$ we have
$$
\mc{G}(H^i_c(X,W(\OO_X))_s)\cong H^i_{\text{\'et},c}(X\times_k \bar{k},\Z_p).
$$
\item For all $i$ and $n\geq 1$ we have
$$
\mc{G}(H^i_c(X,W_n(\OO_X))_s)\cong H^i_{\text{\'et},c}(X\times_k \bar{k},\Z/p^n).
$$
\end{enumerate} 
\begin{proof}
Again, let $Y$ be a compactification of $X$ and $\mathcal{I}$ an ideal for $Z=Y\backslash X$. 
We denote by $\bar{X}:=X\times_k \bar{k}$ the base change, and set $\bar{\mathcal{I}}:=\mathcal{I}\otimes_k \bar{k}$.

Since $W_n(\bar{k})$ is flat over $W_n(k)$  the natural map
\begin{equation} \label{step1-competale}
H^i(Y,W_n(\mathcal{I}))_s \otimes_{W_n(k)} W_n(\bar{k}) \xr{\cong} H^i(\bar{Y},W_n(\bar{\mathcal{I}}))_s
\end{equation}
is an isomorphism.
In view of  Lemma \ref{lemmabasechange}(i) we obtain  
$$
H^i_c(X,\OO_X)_s \otimes_{W(k)} W(\bar{k}) \xr{\cong} H^i_c(\bar{X},\OO_{\bar{X}})_s.
$$

Let $j:\bar{X}\xr{}\bar{Y}$ be the open immersion. 
We have an exact sequence of sheaves on the \'etale site $\bar{Y}_{\text{\'et}}$,
\begin{equation}\label{exseq-competale}
0\xr{} j_! \;\Z/p^n \xr{} W_n(\bar{\mathcal{I}}) \xr{1-F} W_n(\bar{\mathcal{I}}) \xr{} 0.
\end{equation}
Indeed, \ref{exseq-competale} is exact because there is a commutative diagram
$$
\xymatrix
{
0 \ar[r] 
&
i_*\;\Z/p^n \ar[r]
&
W_n(\OO_{\bar{Z}})\ar[r]^{1-F}
&
W_n(\OO_{\bar{Z}})\ar[r]
&
0
\\
0 \ar[r] 
&
\Z/p^n \ar[r] \are[u]
&
W_n(\OO_{\bar{Y}})\ar[r]^{1-F} \are[u]
&
W_n(\OO_{\bar{Y}})\ar[r]\are[u]
&
0 
\\
&
j_! \;\Z/p^n \ar[r] \arir[u]
&
W_n(\bar{\mathcal{I}})\ar[r]^{1-F} \arir[u]
&
W_n(\bar{\mathcal{I}}) \arir[u]
&  
}
$$
with exact first two lines and $i:\bar{Z}=\bar{Y}\backslash \bar{X} \xr{} \bar{Y}$ being the closed immersion.

If $M$ is a finite $W_n(\bar{k})$-module equipped with a $\sigma$-linear endomorphism $F$ then the map $1-F:M\xr{} M$ is surjective \cite[II, Lemme~5.3]{Illusie}. Thus \ref{exseq-competale} yields (2). The compatibility with the Galois action is readily verified. 
The statement (1) follows from the compatibility of $\mc{G}$ with projective limits (see \ref{compGprojlim}).
\end{proof}
\end{proposition}

\begin{remark}
Proposition \ref{proposition-competale}(1) and (2) is well-known for proper schemes. The case $n=1$ is Proposition~2.2.5 in \cite{Katz}. 
\end{remark}

\begin{corollary}
Let $k$ be a perfect field of positive characteristic $p$.
Let $X$ be an \emph{affine} scheme of finite type over $k$.
Suppose $X$ is equidimensional of dimension $d$ and suppose that $X$ is 
Cohen-Macaulay. Then 
$$
H^i_{\text{\'et},c}(X\times_k \bar{k},\Z/p)=0\quad \text{for all $i\neq d$.}
$$
\begin{proof}
This follows from Theorem \ref{thm-weakLefschetz} and Proposition \ref{proposition-competale}.
\end{proof}
\end{corollary}

\section{The formal Euler characteristic of the slope zero part of rigid cohomology modulo powers of $p$.}

\subsection{} 
\assumptionk We denote by $W(k)$ and $K$ the ring of Witt vectors and its quotient field, respectively. 

\begin{notation} \label{notation-lambda} In order to simplify the notation we will write 
$\Lambda$ for $W(k)$ or $W_n(k)$ or $K$. We denote by $\mc{D}_{\Lambda}$ the ring $\mc{D}$ for $\Lambda=W(k)$, the ring $\mc{D}_n$ for $\Lambda=W_n(k)$, 
and the ring $\mc{D}\otimes_{W(k)}K$ for $\Lambda=K$. 
\end{notation}

\begin{definition}\label{definitionKLambda}
 We denote by $\mc{K}_{\Lambda}$ the quotient of 
the free abelian group generated by 
Frobenius stable $\mc{D}_{\Lambda}$-modules $M$ (see Definition \ref{definition-Frobstable}) which are finite and \emph{free} as $\Lambda$-modules, modulo the 
relations
$
[M]-[M']-[M'']
$
for every exact sequence 
$$
0\xr{} M' \xr{} M \xr{} M' \xr{} 0
$$
as $\mc{D}_{\Lambda}$-modules. Note that since $M,M',M''$ are free the sequence is split as sequence of 
$\Lambda$-modules.
\end{definition}
 
In order to get a comparison theorem we need to find natural 
\emph{perfect} complexes  $R\Gamma(X,W\OO_X)_s$ and $R\Gamma(X,W_n\OO_X)_s$ with cohomology groups $H^*_c(X,W\OO_X)_s$ and $H^*_c(X,W_n\OO_X)_s$, respectively. 

\subsection{}
Recall that we have a functor 
$$
(\text{Schemes}/\mathbb{F}_p)\xr{} (\text{Schemes}/\mathbb{F}_p), \quad X\mapsto X^{\perf}, 
$$
defined by $\Spec(A)\mapsto \Spec(A^{\perf})$ on affine schemes, and 
$$
A^{\perf}=\varinjlim_{\Fr} A :=\varinjlim(A\xr{\Fr}A\xr{\Fr}\dots).
$$
In general $X^{\perf}$ in not noetherian, but the underlying topological spaces of $X$ and 
$X^{\perf}$ are identified via the natural map $X^{\perf}\xr{} X$. 

Let $X$ be separated and of finite type over a perfect field $k$. 
Choose a compactification $Y$ of $X$ and an ideal $\mc{I}$ for the complement $Y\backslash X$. 
We define 
$$
H^i_c(X,W_n\OO_{X^{\perf}}) := H^i(Y,W_n(\varinjlim_{\Fr} \mc{I})).
$$
The next proposition implies the independence of $Y,\mc{I}$.

\begin{proposition}\label{proposition-compperfect}
\assumptionk
Let $X$ be separated and of finite type over $k$.
The natural map 
$$
H^i_c(X,W_n\OO_X)_s\xr{} H^i_c(X,W_n\OO_{X^{\perf}}) 
$$
is an isomorphism for all $i$ and all $n\geq 1$.
\begin{proof}
Since the topological space of $Y^{\perf}$ is noetherian we obtain
$$
H^i(Y,W_n(\varinjlim_{\Fr} \mc{I}))= H^i(Y,\varinjlim_{F} W_n(\mc{I})) =\varinjlim_F H^i(Y,W_n(\mc{I})).
$$ 
In view of Lemma \ref{lemma-stable/nil} we see that
$$
H^i(Y,W_n(\mc{I})) \xr{} \varinjlim_F H^i(Y,W_n(\mc{I}))
$$
induces an isomorphism with the stable part $H^i(Y,W_n(\mc{I}))_s$. 
\end{proof} 
\end{proposition}

\subsection{} An $\F_p$-algebra $R$ is called perfect if $\Fr$ induces an automorphism
of $R$, equivalently $R=R^{\perf}$; an ideal $I\subset R$ is called perfect if $\Fr$ induces
an automorphism on $I$. For any $\F_p$-algebra $R$ and any ideal $I\subset R$ the ideal 
$
\varinjlim_{\Fr} I\subset R^{\perf}
$ 
is perfect.

\begin{lemma}\label{lemma-reductionperfect}
Let $R$ be a perfect $k$-algebra. Let $I\subset R$ be a perfect ideal.
The following holds:
\begin{itemize}
\item[(i)] For all $m\geq 1$: $W(I)\otimes_{W(k)}W_m(k)=W_m(I)$.
\item[(i.n)] For all $n\geq m\geq 1$: $W_n(I)\otimes_{W_n(k)}W_m(k)=W_m(I)$.
\item[(ii)] The $W(k)$-module $W(I)$ is flat.
\item[(ii.n)] The $W_n(k)$-module $W_n(I)$ is flat for all $n\geq 1$.
\end{itemize}
\begin{proof}
We leave the proof to the reader.
\end{proof}
\end{lemma}

\subsection{}
Let $X$ be separated and of finite type over a perfect field $k$. 
Choose a compactification $Y$ of $X$ and an ideal $\mc{I}$ for the complement $Y\backslash X$. Attached to a finite affine covering $\{U_i\}_i$ of $Y$
we get the \v{C}ech-complex 
\begin{equation}\label{Rgamma}
R\Gamma(\{U_i\},W(\varinjlim_{\Fr} \mc{I}))\in D^b(\text{$W(k)$-modules}).
\end{equation}
We have 
$$
R\Gamma(\{U_i\},W(\varinjlim_{\Fr} \mc{I})) = \varprojlim_n R\Gamma(\{U_i\},W_n(\varinjlim_{\Fr} \mc{I})) = R\varprojlim_n R\Gamma(\{U_i\},W_n(\varinjlim_{\Fr} \mc{I})).
$$
Indeed, the first equality is obvious and the second follows because 
$$
R\Gamma(\{U_i\},W_n(\varinjlim_{\Fr} \mc{I}))\xr{} R\Gamma(\{U_i\},W_m(\varinjlim_{\Fr} \mc{I}))
$$
is surjective on the components.
Proposition \ref{proposition-compperfect} implies
\begin{multline*}
H^i(R\Gamma(\{U_i\},W(\varinjlim_{\Fr} \mc{I})))=H^iR\varprojlim_n R\Gamma(\{U_i\},W_n(\varinjlim_{\Fr} \mc{I}))\\
\overset{(*)}{=} \varprojlim_n H^i_c(X,W_n\OO_X)_s = H^i_c(X,W\OO_X)_s
\end{multline*}
for all $i$. For $(*)$ we used 
$$
H^iR\Gamma(\{U_i\},W_n(\varinjlim_{\Fr} \mc{I}))=H^i(Y,W_n(\varinjlim_{\Fr} \mc{I})),
$$
and the fact that the projective system 
$$
n\mapsto H^i(Y,W_n(\varinjlim_{\Fr} \mc{I}))\cong H^i_c(X,W_n\OO_X)_s
$$
satisfies the Mittag-Leffler condition. 

Thus the complex \ref{Rgamma} is in $D^b(\text{$W(k)$-modules})$ up to 
canonical isomorphisms independent of the choice of $Y,\mc{I},$ and the covering $\{U_i\}$.
\begin{definition} We define 
$$
R\Gamma_c(X,W\OO_X)_s:=R\Gamma(\{U_i\},W(\varinjlim_{\Fr} \mc{I})),
$$
and similarly for all $n\geq 1$,
$$
R\Gamma_c(X,W_n\OO_X)_s:=R\Gamma(\{U_i\},W_n(\varinjlim_{\Fr} \mc{I})).
$$
If $X$ is proper we drop the index $c$.
\end{definition}

Since the components of \ref{Rgamma} are flat by Lemma \ref{lemma-reductionperfect}, 
$R\Gamma_c(X,W\OO_X)$ is a perfect complex (see \cite[I]{IllusieSGA6I} for the definition of 
a perfect complex which should not be confused with the notion of a perfect ring) and 
\begin{equation}\label{RGammareduction}
R\Gamma_c(X,W\OO_X)_s\otimes^{\mathbb{L}}_{W(k)} W_n(k)=R\Gamma_c(X,W_n\OO_X)_s. 
\end{equation}
The Frobenius morphism 
$$
F:R\Gamma(\{U_i\},W(\varinjlim_{\Fr} \mc{I}))\xr{} \sigma_*R\Gamma(\{U_i\},W(\varinjlim_{\Fr} \mc{I}))
$$
induces an isomorphism in $D^b(\text{$W(k)$-modules})$:
$$
F:R\Gamma_c(X,W\OO_X)_s \xr{} \sigma_*R\Gamma_c(X,W\OO_X)_s,
$$
which agrees with our $F$-operation on the cohomology. 

\subsection{}
As a unifying notation we write $\Lambda$ for $W(k)$ and $W_n(k)$ in the following.
Suppose $N\in D^b(\text{$\Lambda$-modules})$ is a perfect complex together 
with a quasi-isomorphism $F:N\xr{} \sigma_*N$. By definition we can find 
$$M\in K^b(\text{free and finite $\Lambda$-modules})=:K^b(\text{ff-$\Lambda$})$$ together with 
$F_M:M\xr{} \sigma_*M$ and a quasi-isomorphism 
$\psi_M:M\xr{} N$ such that $\sigma_*(\psi_M)\circ F_M=F\circ \psi_M$ in $D^b(\Lambda)$. 
Induced by $F_M$ we get a $\sigma$-linear map $F_M^i$ on the components $M^i$ of $M$. 
Taking the Frobenius-stable part is an exact functor when restricted to finite $\Lambda$-modules
and preserves free modules (see Proposition \ref{proposition-stablexact}), therefore $M_s\xr{} M$
is an isomorphism in $K^b(\text{ff-$\Lambda$})$. For all $i$ we get a stable $\mc{D}_{\Lambda}$-module 
$M^i_s$ with $F$-operation induced by $F^i_M$, and $V=pF^{-1}$. 
We define
\begin{equation}\label{eulerclass}
e(M,F_M)=\sum_i (-1)^i [M^i_s]\in \mc{K}_{\Lambda}.
\end{equation}

\begin{proposition}
The class $e(M,F_M,\psi_M)$ depends only on $N$ and the morphism $F:N\xr{} \sigma_*N$.
\begin{proof}
First, we observe that if $M\cong 0$ in $K^b(\text{ff-$\Lambda$})$ then $e(M,F_M)=0$
for any $F_M$.

Let $M_1,M_2\in K^b(\text{ff-$\Lambda$})$ and $F_{M_i}:M_i\xr{} \sigma_*M_i, i=1,2,$
with a morphism of complexes $\psi:M_1\xr{} M_2$ such that $\psi$ is an isomorphism in $K^b(\text{ff-$\Lambda$})$ and $\sigma_*(\psi)\circ F_{M_1}=F_{M_2}\circ \psi$ in 
$K^b(\text{ff-$\Lambda$})$. In particular, there is a homotopy 
$K:M_1\xr{} \sigma_*M_2[-1]$ such that 
$$
F_{M_2}\circ \psi-\sigma_*\psi\circ F_{M_1}=K\circ d_{M_1}+d_{M_2}\circ K.
$$
We need to show that 
$
e(M_1,F_{M_1})=e(M_2,F_{M_2}). 
$ 

We define a morphism of complexes 
$$
F_3:{\rm cone}(\psi)\xr{} {\rm cone}(\sigma_*\psi)=\sigma_*{\rm cone}(\psi)
$$
where $F^i_3:M_2^i\oplus M_1^{i+1}\xr{} \sigma_*M_2^i\oplus \sigma_*M_1^{i+1}$ is
of the form 
$$
F^i_3=\left(\begin{matrix} F^i_2 & K^{i+1} \\ 0 & F^{i+1}_1\end{matrix}\right).
$$
Therefore we get short exact sequences which are compatible with the $F$-operations
\begin{equation*}
\begin{split}
&0\xr{} (M_2^i,F^i_2)\xr{} (M_2^i\oplus M_1^{i+1}, F_3) \xr{}  (M_1^{i+1},F^{i+1}_1) \xr{} 0 \\
&0\xr{} (M_2^i,F^i_2)_s\xr{} (M_2^i\oplus M_1^{i+1}, F_3)_s \xr{}  (M_1^{i+1},F^{i+1}_1)_s \xr{} 0. 
\end{split}
\end{equation*}
Since ${\rm cone}(\psi)\cong 0$ this implies the claim.
\end{proof}
\end{proposition} 

In view of this proposition we write $e(N,F)$ for the class constructed in \ref{eulerclass}.

\subsection{} 
Recall that we have obvious maps:
$$
\xymatrix
{
\mc{K}_{W(k)} \ar[r]^{\otimes_{W(k)}K} \ar[d]^{\otimes_{W(k)}W_n(k)}
&
\mc{K}_K
\\
\mc{K}_{W_n(k)}
}
$$

\begin{proposition}\label{proposition-ecomp}
Let $X$ be separated and of finite type over $k$.
\begin{itemize}
\item[(i)] For all $n\geq 1$:
$$
e(R\Gamma_c(X,W\OO_X)_s,F)\otimes_{W(k)}W_n(k)=e(R\Gamma_c(X,W_n\OO_X)_s,F).
$$ 
\item[(ii)] The following equality holds
$$
e(R\Gamma_c(X,W\OO_X)_s,F)\otimes_{W(k)}K=\sum_i (-1)^i [(H^i_c(X,W\OO_X)_s\otimes_{W(k)}K]\in \mc{K}_{K}.
$$ 
\item[(iii)] Let $n\geq 1$. If the cohomology groups $H^i_c(X,W_n\OO_X)_s$ are free 
$W_n(k)$-modules for all $i$ then 
$$
e(R\Gamma_c(X,W_n\OO_X)_s,F)=\sum_i (-1)^i [(H^i_c(X,W_n\OO_X)_s]\in \mc{K}_{W_n(k)}.
$$
\end{itemize}
\begin{proof}
Statement (i) follows from \ref{RGammareduction}.  Statements (ii) and (iii) are obvious.
\end{proof}
\end{proposition}



\subsection{}
Let $k$ be a finite field with $q=p^a$ elements. Let $r$ be an integer $r\geq 1$.
We define a homomorphism 
$$
{\rm Tr}^r_{\Lambda}:\mc{K}_{\Lambda}\xr{} \Lambda, \quad [M]\mapsto {\rm Tr}(F^{ar}\mid M).
$$
Obviously, for an element $M\in \mc{K}_{W(k)}$ we get 
$$
{\rm Tr}^r_K(M\otimes_{W(k)} K)={\rm Tr}^r_{W(k)}(M), \quad {\rm Tr}^r_{W_n(k)}(M\otimes_{W(k)} W_n(k))\equiv {\rm Tr}^r_{W(k)}(M) \mod p^n.
$$

\begin{corollary}\label{corollary-congruenceformulaformal}
Let $k$ be a finite field with $q=p^a$ elements. Let $X$ be a proper scheme over $k$. We denote by $H^*_{rig}(X/K)_{[0]}$ the slope zero part of rigid cohomology.
For an integer $r\geq 1$ we denote by $k_r$ the field extension of $k$ of degree $r$. 
Let $n\geq 1$ be an integer.  
\begin{itemize}
\item[(i)] The trace
$$
N:=\sum_{i\geq 0} (-1)^i {\rm Tr}(F^{ar}\mid H^i_{rig}(X/K)_{[0]})
$$
lies in $W(k)$, i.e. $N\in W(k)$.
For all $r,n$ the following congruence holds:
$$
N\equiv {\rm Tr}^r_{W_n(k)}(e(R\Gamma(X,W_n\OO_X)_s,F)) \mod p^{n}.
$$ 
\item[(ii)] There is an integer $r_0$ which depends on $n$ and $X$ such that for all $r\geq r_0$
the following congruence holds:
$$
{\rm Tr}^r_{W_n(k)}(e(R\Gamma(X,W_n\OO_{X})_s,F)) \equiv \#X(k_r) \mod p^{n}. 
$$
\item[(iii)] Suppose that $H^i(X,\OO_X)=H^i(X,\OO_X)_s$ for all $i$. 
Then 
$$
{\rm Tr}^r_{W_n(k)}(e(R\Gamma(X,W_n\OO_X)_s,F)) \equiv \#X(k_r) \mod p^{\min\{ra,n\}}.
$$ 
\end{itemize}
\begin{proof}
For (i): Proposition \ref{proposition-ecomp} implies 
\begin{equation*}
\begin{split}
{\rm Tr}^r_{W(k)}(e(R\Gamma(X,W\OO_X)_s,F)) &= \sum_{i\geq 0} (-1)^i {\rm Tr}(F^{ar}\mid H^i(X,W\OO_X)_s\otimes_{W(k)}K), \\
{\rm Tr}^r_{W(k)}(e(R\Gamma(X,W\OO_X)_s,F)) &\equiv {\rm Tr}^r_{W_n(k)}(e(R\Gamma(X,W_n\OO_X)_s,F)) \mod p^{n}. 
\end{split}
\end{equation*}
Thus the claim follows from Proposition \ref{proposition-slopezero}.

For (ii): We have the trace formula for rigid cohomology:
$$
\sum_{i\geq 0} (-1)^i {\rm Tr}(F^{ar}\mid H^i_{rig}(X/K))=X(k_r).
$$ 
The isocrystal $H^*_{rig}(X/K)$ has a slope decomposition with slopes $\geq 0$. 
On the slope $\lambda$-part the eigenvalues of $F^a$ have $p$-adic valuation 
$a\cdot \lambda$.
Let $\lambda>0$ be the smallest slope $\neq 0$. Choose $r_0$ such that $a\cdot r_0\cdot \lambda>n$. For every $r\geq r_0$
we get 
$$
N:=\sum_{i\geq 0} (-1)^i {\rm Tr}(F^{ar}\mid H^i_{rig}(X/K)_{[0]})=X(k_r) \mod  p^n,
$$
where $N\in W(k)$ by part (i). In view of (i), this proves the claim.

For (iii): Proposition \ref{proposition-FstableOimpliesFstableWO} implies that 
$$
H^i(X,W\OO_X)_s=H^i(X,W\OO_X) \quad \text{for all $i$.}
$$
In view of Proposition \ref{proposition-ecomp}(ii) we see that 
$$
N:= {\rm Tr}^r_{W(k)}(e(R\Gamma(X,W\OO_X)_s,F)) = \sum_{i\geq 0} (-1)^i {\rm Tr}(F^{ra}\mid H^i(X,W\OO_X)\otimes_{W(k)}K) 
$$
lies in $W(k)$. It is proved in \cite[Corollary~1.3]{BBE} that $N\equiv \#X(k_r) \mod p^{ar}$. Now,  Proposition \ref{proposition-ecomp}(i) implies
$$
N \equiv {\rm Tr}^r_{W_n(k)}(e(R\Gamma(X,W_n\OO_X)_s,F)) \mod p^n,
$$
which proves the statement.
\end{proof}
\end{corollary}

\begin{corollary}\label{corollary-congruenceformula}
Let $k$ be a finite field with $q=p^a$ elements. Let $X$ be a proper scheme over $k$. 
Let $n\geq 1$ be an integer and suppose that $H^i(X,W_n\OO_X)$ is a free $W_n(k)$-module
for all $i$. Then 
$$
\sum_{i\geq 0} (-1)^i {\rm Tr}(F^a\mid H^i(X,W_n\OO_X)) \equiv \#X(k) \mod p^{\min\{a,n\}}.
$$ 
\begin{proof}
The case $n=1$ is Katz's formula \cite{Katz}. In the case $n\geq 2$, 
Proposition \ref{proposition-properfree} 
implies  
$
H^i(X,\OO_X)_s=H^i(X,\OO_X) 
$
for all $i$.
Proposition \ref{proposition-ecomp}(iii) and Corollary \ref{corollary-congruenceformulaformal}(iii)
imply the claim.
\end{proof}
\end{corollary}

\begin{example}
Let $k$ be finite field with $p^a$ elements. Let $G$ be a finite (abstract) group. For every integer $d\geq 1$, Serre shows the existence of a regular complete 
intersection $Y$ in projective space such that $\dim(Y)=d$ and $G$ acts freely on $Y$ \cite[\textsection20]{Serre}.
 
Since $Y$ is a complete intersection, we get for all $n\geq 1$: 
\begin{equation*}
H^i(Y,W_n\OO_Y)=0 \quad \text{for all $i\not\in \{0,d\}$}, \quad H^0(Y,W_n\OO_Y)=W_n(k).
\end{equation*}
Because $W_n\OO$ is an \'etale sheaf and \'etale cohomology for $W_n\OO$ agrees with Zariski cohomology \cite[Proposition~0.1.5.8]{Illusie},
we obtain a spectral sequence 
\begin{equation}\label{equation-spsec}
E_2^{p,q}=H^p(G,H^q(Y,W_n\OO_Y))\Rightarrow H^{p+q}(X,W_n\OO_X),
\end{equation}
where $X:=Y/G$.
From the spectral sequence we obtain
\begin{equation}\label{HiHiG}
H^i(X,W_n\OO_X)=H^i(G,W_n(k))\quad \text{for all $i<d$.}
\end{equation}
In particular, $H^i(X,W_n\OO_X)=H^i(X,W_n\OO_X)_s$ for all $i<d$. 

For $G=\Z/p^m$ we have 
\begin{equation}\label{HiZ/pm}
H^i(G,W_n(k))=\begin{cases} W_n(k) &\text{for $i=0$,}\\ 
\ker(p^m:W_n(k)\xr{} W_n(k)) &\text{for $i$ odd,} \\
\coker(p^m:W_n(k)\xr{} W_n(k)) &\text{for $i$ even.}
\end{cases}
\end{equation}

Suppose $d\geq 2$. We conclude that the cohomology groups $H^*(X,W_n\OO_X)_s$ are free $W_n(k)$-modules 
if and only if $m\geq n$ (for $H^d$ we use Proposition \ref{proposition-equivalencefreestable}). For the projective limit we have the following picture: 
\begin{equation*}
H^i(X,W\OO_X)_s=\begin{cases} 0 &\text{if $i$ is odd and $i<d$,} \\ W_m(k) &\text{if $i$ is even and $0<i<d$.}\end{cases}
\end{equation*}
For $i=d$ we have to distinguish two cases:
\begin{enumerate}
\item If $d$ is odd then $H^d(X,W\OO_X)_s$ is free of rank $=\dim H^d(X,\OO_X)_s$.
\item If $d$ is even then  $H^d(X,W\OO_X)_s$ has $W_m(k)$ as torsion subgroup. The dimension of $H^d(X,W\OO_X)_s\otimes \Q$ is $\dim H^d(X,\OO_X)_s-1$. 
\end{enumerate}
This follows from the long exact sequence in Proposition \ref{proposition-longexseq}.
Suppose $d$ is even. We claim that for $r$ sufficiently large 
and $m\geq n$:
\begin{equation*}
\begin{split}
{\rm Tr}(F^{ra}\mid H^d(X,W\OO_X)\otimes \Q)&=
{\rm Tr}(F^{ra}\mid H^d(X,W\OO_X)_s\otimes \Q)\\
&\equiv {\rm Tr}(F^{ra}\mid H^d(X,W_n\OO_X)_s)-1 \mod p^n.
\end{split}
\end{equation*}
The first equality follows (for $r\gg0$) because $ H^d(X,W\OO_X)_s\otimes \Q$ is the
slope zero part of $H^d(X,W\OO_X)\otimes \Q$ (see the proof of Corollary \ref{corollary-congruenceformulaformal}). The second equality is a consequence of 
Proposition \ref{proposition-ecomp}:
\begin{multline*}
1+ {\rm Tr}(F^{ra}\mid H^d(X,W\OO_X)_s\otimes \Q) =
\sum_{i\geq 0}(-1)^i  {\rm Tr}(F^{ra} \mid H^i(X,W\OO_X)_s\otimes \Q)\equiv\\
\sum_{i\geq 0}(-1)^i  {\rm Tr}(F^{ra} \mid H^i(X,W_n\OO_X)_s) 
={\rm Tr}(F^{ra}\mid H^d(X,W_n\OO_X)_s) \mod p^n.
\end{multline*}
For the last equality we used \ref{HiHiG}, \ref{HiZ/pm}, and the assumption $m\geq n$. 
\end{example}


\end{document}